\documentclass[11pt]{amsart} 
\usepackage{amssymb,amsmath,latexsym} 

\newtheorem{theorem}{Theorem} 
\newtheorem{lemma}[theorem]{Lemma} 
\newtheorem{ex}[theorem]{Example}

\title{On Stanley's reciprocity theorem for rational cones}

\author{Matthias Beck}
\address{Department of Mathematics, San Francisco State University, 
San Francisco, CA 94132, U.S.A.}
\email{beck@math.sfsu.edu}
\urladdr{http://math.sfsu.edu/beck}

\author{Mike Develin}
\address{American Institute of Mathematics, Palo Alto, CA 94306, U.S.A.}
\email{develin@post.harvard.edu}

\begin{document}
\setlength{\parindent}{0pt} 

\def\K{\mathcal{K}}
\def\P{\mathcal{P}}
\def\Q{\mathbb{Q}}
\def\R{\mathbb{R}}
\def\Z{\mathbb{Z}}
\def\SL{\mathit{SL}}
\def\A{{\bf A}} 
\def\B{{\bf B}} 
\def\a{{\bf a}}
\def\b{{\bf b}}
\def\m{{\bf m}}
\def\r{{\bf r}}
\def\v{{\bf v}}
\def\w{{\bf w}}
\def\x{{\bf x}}
\def\y{{\bf y}}
\def\z{{\bf z}}
\newcommand\const{\operatorname{const}} 

\abstract 

We give a short, self-contained proof of Stanley's reciprocity theorem for a 
rational cone $\K \subset \R^d$. Namely, let $\sigma_\K (\x) = \sum_{ \m \in \K \cap \Z^d } \x^\m$. 
Then $\sigma_\K (\x)$ and $\sigma_{ \K^\circ } (\x)$ are rational functions which satisfy the 
identity $\sigma_\K (1/\x) = (-1)^d \sigma_{ \K^\circ } (\x)$. A corollary of Stanley's theorem is 
the Ehrhart-Macdonald reciprocity theorem for the lattice-point enumerator of rational polytopes. A 
distinguishing feature of our proof is that it uses neither the shelling of a polyhedron nor the 
concept of finite additive measures.
The proof follows from elementary techniques 
in contour integration.

\endabstract 


\thanks{Thanks to Sinai Robins and Guoce Xin for helpful discussions.}
\date{\today}

\maketitle  

\setlength{\parskip}{0.4cm}
\bibliographystyle{amsplain}


\section{Introduction}

Let $\K$ be a \emph{rational cone}, that is, of the form
\[
  \K = \left\{\r \in \R^d : \ \A \, \r \le 0 \right\}
\]

for some integral $m \times d$-matrix $\A$. Thus $\K$ is a real cone in 
$\R^d$ defined by $m$ inequalities.
 We assume that $\K$ is $d$-dimensional and pointed: i.e., that $\K$ does not 
contain a line.

In his study of nonnegative integral solutions to linear systems, Stanley was led to consider the generating functions
\[
  \sigma_\K (\x) = \sum_{ \m \in \K \cap \Z^d } \x^\m
\]
and its companion $\sigma_{ \K^\circ } (\x)$ for the interior $\K^\circ$ of $\K$. Here $\x^\m$ denotes the product $x_1^{ m_1 } x_2^{ m_2 } \cdots x_d^{ m_d }$. The function $\sigma_\K$ (as well as $\sigma_{ \K^\circ }$) is a rational function in the components of $\x$; this is proved, for example, by triangulating $\K$ into simplicial cones, for which one can explicitly form the rational functions representing their generating functions (see, for example, \cite[Chapter 4]{stanleyec1}).
The fundamental reciprocity theorem of Stanley \cite{stanleyreciprocity} relates the two rational functions $\sigma_\K$ and $\sigma_{ \K^\circ }$. We abbreviate the vector $\left( 1/x_1, 1/x_2, \dots, 1/x_d \right)$ by $1/\x$.

\begin{theorem}[Stanley]\label{stanrec} 
As rational functions, $\sigma_\K \left( \x \right) = (-1)^d \, \sigma_{ \K^\circ } (1/\x)$.
\end{theorem}

The proof of Stanley's theorem is not quite as simple as that of the rationality of $\sigma_\K$. Most 
proofs of Theorem \ref{stanrec} that we are aware use either the concept of shelling of a polyhedron 
or the concept of finitely additive measures (also called valuations). Both concepts reduce the 
theorem to simplicial cones, for which Theorem \ref{stanrec} is not hard to prove. In this paper, we 
give a proof which neither depends on shelling nor on valuations. It is based on rational generating 
functions, which essentially go back to Euler. In \cite{stanleyreciprocity} Stanley used these 
generating functions to suggest a residue-computation approach to proving Theorem \ref{stanrec}; he 
then implemented it by computing each of the residues in gory detail. Here we show that one can, in 
fact, reduce Theorem \ref{stanrec} to a simple change of variables in the integral over the 
generating function, proving Stanley's reciprocity theorem without any actual computation of 
residues or other technical machinery. 


\section{Euler's generating function}

Suppose we have an element $\phi$ of $\SL_d(\Z)$, i.e. a unimodular transformation preserving the integer lattice. We 
claim that the image $\phi(\K)$ will satisfy the reciprocity law if and only if $\K$ does. Indeed, it is easy to see 
that the generating function $\sigma_{\phi(\K)}$ is obtained by substituting $\x^{\phi(e_i)}$ for $x_i$ (for all 
$i$) in the generating function $\sigma_{\K}$, and similarly the generating function $\sigma_{\phi(\K^\circ)}$ is 
obtained from $\sigma_{\K^\circ}$ by the same substitutions. From this observation, we immediately obtain that if 
$\sigma_{\K}(\x)=\sigma_{\K^\circ}(1/\x)$, then
\[
\sigma_{\phi(\K)}(\x) = \sigma_{\K} \left( \x^{\phi(e_i)} \right) = \sigma_{\K^\circ} \left( \x^{-\phi(e_i)} \right) = 
\sigma_{\phi(\K^\circ)}(1/\x),
\]
so $\sigma_{\phi(\K)}$ and $\sigma_{\phi(\K^\circ)}$ also satisfy the reciprocity law.

Now, given any pointed cone $\K$, we can find a unimodular transformation $\phi$ so that $\phi(\K)$ lies in the 
nonnegative orthant, and intersects its boundary only at $(0,0,\ldots,0)$: simply pick a lattice basis such that 
$\K$ is contained in the interior of the cone it spans, and send this to the standard basis. (This is easy to do; 
for instance, take any hyperplane which intersects $\K$ only at the origin, pick a lattice basis for 
the sublattice 
contained in it, and take the final point at lattice distance 1 in the direction of $\K$, and very 
far away in the 
direction opposite the $d-1$ points in the hyperplane.) Therefore, it suffices to consider the case where $\K\subset 
\R_{\ge 0}^d$ and $(0,0,\ldots,0)$ is the only point in $\K$ with any zero coordinate. We will assume that $\K$ has 
this form from now on.

Denote the 
columns of $\A$ by $\a_1, \dots, \a_d$. Our main tool for proving Theorem 
\ref{stanrec} is the following lemma, the idea of which goes back to at 
least Euler \cite{euler} and which is proved by simply expanding geometric 
series.

\begin{lemma}\label{eulerlemma}
$\sigma_\K (\x)$ equals the constant $\z$-coefficient of the function 
\[ 
  \frac 1 { \left( 1 - x_1 \z^{ \a_1 } \right) \left( 1 - x_2 \z^{ \a_2 } \right) \cdots \left( 1 - x_d \z^{ \a_d } \right) \left( 1 - z_1 \right) \left( 1 - z_2 \right) \cdots \left( 1 - z_m \right) } 
\] 
expanded as a power series centered at $\z = 0$. 
\end{lemma} 

We will use an integral version of this lemma, for which we need additional variables to avoid integrating over singularities. Let
\begin{equation}\label{xyGeneratingFunction}  \theta_\K (\x,\y) = \theta_\K \left( x_1, \dots, x_d, y_1, \dots, y_m \right) = \const_\z \left( \frac 1 { \prod_{ j=1 }^{ d } \left( 1 - x_j \z^{ \a_j } \right) \prod_{ k=1 }^{ m } \left( 1 - y_k z_k \right) } \right) .
\end{equation}
This is a rational function in the coordinates of $\x$ and $\y$, and 
$\sigma_\K (\x) = \theta_\K (\x, 1, 1, \dots, 1)$. The $y$-variables 
represent slack variables, one for each inequality; thus instead of $\A\,\x \le 0$ and $\x\ge 0$, we 
consider the cone given by $(\A\,|\,I)\begin{pmatrix} \x \\ \y \end{pmatrix} = 0$ and $\x, \y \ge 0$, 
where 
$I$ is the $n$-by-$n$ identity matrix. We next translate 
this 
constant-term identity into an integral identity. We abbreviate the 
measure $\frac{ dz_1 }{ z_1 } \frac{ dz_2 }{ z_2 } \cdots \frac{ dz_m }{ 
z_m }$ by $\frac{ d\z }{ \z }$, and we write $|\x|<1$ to indicate that 
$\left| x_1 \right| , \left| x_2 \right| , \dots , \left| x_d \right| < 
1$. Then for $|\x|, |\y| < 1$
\[
  \theta_\K (\x,\y) = \int \frac 1 { \prod_{ j=1 }^{ d } \left( 1 - x_j \z^{ \a_j } \right) \prod_{ 
k=1 }^{ m } \left( 1 - y_k z_k \right) } \frac{ d \z }{ \z } \ ,
\]
where the integral sign stands for $1/(2 \pi i)^m$ times an $m$-fold 
integral, each one over the unit circle. Note that although we are going 
to plug in $\y=(1,1,\ldots,1)$, we may compute the rational function 
$\theta_\K(\x,\y)$ by evaluating it in any sufficiently large domain, 
such as the one used here. 

Analogously, we have for the open generating function
\[
  \sigma_{ \K^\circ } (\x) = \const \left( \frac{ x_1 \z^{ \a_1 } }{ 1 - x_1 \z^{ \a_1 } } \frac{ x_2 \z^{ \a_2 } }{ 1 - x_2 \z^{ \a_2 } } \cdots \frac{ x_d \z^{ \a_d } }{ 1 - x_d \z^{ \a_d } } \frac{ z_1 }{ 1-z_1 } \frac{ z_2 }{ 1-z_2 } \cdots \frac{ z_m }{ 1-z_m } \right) 
\]
and we define
\[
  \theta_{\K^\circ} (\x,\y) = \const_\z \left( \prod_{ j=1 }^{ d } \frac{ x_j \z^{ \a_j } }{ 1 - x_j \z^{ \a_j } } \prod_{ k=1 }^{ m } \frac{ y_k z_k }{ 1 - y_k z_k } \right) .
\]
Then $\sigma_{\K^\circ} (\x) = \theta_{\K^\circ} (\x, 1, 1, \dots, 1)$ and for $|\x|, |\y| < 1$
\[
  \theta_{\K^\circ} (\x,\y) = \int \prod_{ j=1 }^{ d } \frac{ x_j \z^{ \a_j } }{ 1 - x_j \z^{ \a_j } } \prod_{ k=1 }^{ m } \frac{ y_k z_k }{ 1 - y_k z_k } \frac{ d \z }{ \z } \ .
\]
The integral representations of $\theta_\K$ and $\theta_{ \K^\circ }$ now suggest how to prove Theorem \ref{stanrec}---make a change of variables $\z \to 1/\z$, say in $\theta_{ \K^\circ }$:
\[
  \theta_{\K^\circ} (\x,\y) = \int \prod_{ j=1 }^{ d } \frac{ x_j \z^{ -\a_j } }{ 1 - x_j \z^{ -\a_j } } \prod_{ k=1 }^{ m } \frac{ y_k z_k^{ -1 } }{ 1 - y_k z_k^{ -1 } } \frac{ d \z }{ \z } \ .
\]
Hence the rational function $\theta_{ \K^\circ } (1/\x,1/\y)$ has the integral representation
\begin{align*}
  \theta_{\K^\circ} (1/\x,1/\y) &= \int \prod_{ j=1 }^{ d } \frac{ x_j^{-1} \z^{ -\a_j } }{ 1 - x_j^{-1} \z^{ -\a_j } } \prod_{ k=1 }^{ m } \frac{ y_k^{-1} z_k^{ -1 } }{ 1 - y_k^{-1} z_k^{ -1 } } \frac{ d \z }{ \z } \\
                                &= \int \prod_{ j=1 }^{ d } \frac{ 1 }{ x_j \z^{ \a_j } - 1 } \prod_{ k=1 }^{ m } \frac{ 1 }{ y_k z_k - 1 } \frac{ d \z }{ \z } \\
                                &= (-1)^{ d+m } \int \prod_{ j=1 }^{ d } \frac{ 1 }{ 1 - x_j \z^{ \a_j } } \prod_{ k=1 }^{ m } \frac{ 1 }{ 1 - y_k z_k } \frac{ d \z }{ \z } \ ,
\end{align*}
valid for $|\x|, |\y| > 1$.
It remains to prove that the rational function given by the integral
\begin{equation}\label{firstint}  \int \frac 1 { \prod_{ j=1 }^{ d } \left( 1 - x_j \z^{ \a_j } \right) \prod_{ k=1 }^{ m } \left( 1 - y_k z_k \right) } \frac{ d \z }{ \z }
\end{equation}with $|\x|, |\y| < 1$ equals the rational function given by the integral
\begin{equation}\label{secondint}  (-1)^{ m } \int \frac 1 { \prod_{ j=1 }^{ d } \left( 1 - x_j \z^{ \a_j } \right) \prod_{ k=1 }^{ m } \left( 1 - y_k z_k \right) } \frac{ d \z }{ \z }
\end{equation}with $|\x|, |\y| > 1$, since once we have this equality, we 
have $\theta_{\K^\circ}(1/\x,1/\y) = (-1)^d \theta_{\K}(\x,\y)$; 
plugging in $\y=(1,1,\ldots,1)$ (whereupon $1/\y = (1,1, \ldots, 1)$ as 
well) shows that $\theta_{\K^\circ}(1/\x) = \theta_{\K}(\x)$ as desired.

Let us consider as the innermost integral, the one with respect to $z_1$. Almost all
of the poles of the integrand $f(z_1)$ are at the solutions $z_1$ of the equations $1 - x_j
\z^{ \a_j } = 0$ (for $j=1, \dots, d$) and $1 - y_1 z_1 = 0$. Since $\left| z_2
\right| , \left| z_3 \right| , \dots, \left| z_m \right| = 1$, each $z_1$-pole is
inside or outside the unit circle, depending on the exponent of $z_1$ in $\z^{ \a_j
}$ and on whether a given $x_j$ or $y_j$ has magnitude smaller or larger than 1.
But this means that the $z_1$-integrals in \eqref{firstint} and \eqref{secondint}
pick up the residues of complementary singularities.

The only other potential poles are at zero and infinity, induced by the 
extra factor of $\frac{1}{z_1}$. We claim that there are no residues at these poles, if they exist. Indeed,
as $z_1$ approaches zero, the residue $z_1 f(z_1)$ is the product of the factors $\frac{1}{1-x_j \z^{\a_j}}$ and 
$\frac{1}{1-y_k z_k}$. Since $K$ is in the nonnegative orthant, in each inequality at least one coefficient 
must be negative; therefore, in at least one of the former factors, the exponent of $z_1$ must be negative. This 
factor then goes to zero as $z_1$ does, while the norms of all of the other factors either go to zero (if the 
exponent of $z_1$ is negative), one (if the exponent is positive), or a constant (if $z_1$ does not appear at all.) 
Therefore, this residue is equal to zero.

A similar argument eliminates the residue at infinity. After a change of variables $z_1\rightarrow
\frac{1}{z_1}$, this residue at infinity is (up to sign) the limit of the product of these same
factors. However, the factor $\frac{1}{1-y_1 z_1}$ goes to zero, while all other factors again go to
zero, one, or a constant, depending on whether the exponent of $z_1$ is positive, negative, or zero
respectively. This completes the argument that the $z_1$-integrals in \eqref{firstint} and
\eqref{secondint} differ by a minus sign.

We can use the same argument for the next variable if we know that the 
$z_1$-integral results in a rational function with a similar-looking 
denominator as the integrands in \eqref{firstint} and \eqref{secondint}. 
But this follows from Euler's generating function: The $z_1$-integral in 
\eqref{firstint} gives the generating function of the cone described by 
$\r \ge 0$ and the first row inequality in $\A \, \r \le 0$, in the 
variables $x_1 \left( z_2, \dots, z_m \right)^{ \a_1' }, \dots, x_d \left( 
z_2, \dots, z_m \right)^{ \a_d' }$ and $\y$, where $\a_1', \dots, \a_d'$ 
are the column vectors of $\A$ after we removed the first row. It is not 
hard to show (from the simplicial case) that the generating function of 
any such cone has as denominator a product of terms of the form 1 minus a 
monomial of the variables, so the $z_2$-integrand has the desired form.

There is a possible obstruction here, which is that these monomials 
may have $x$'s and $y$'s appearing both in the numerator and the 
denominator. If we pick a poor choice of norms of the $x$'s and $y$'s, we 
may not be able to evaluate these integrals, since the norm of this 
monomial may be 1, in which case there are poles on the unit circle. 
However, we can get around this simply by looking at an open set where 
none of these magnitudes are 1 to evaluate the integral for $|\x|, |\y| > 
1$, and looking at the image of this set under $\x\rightarrow \x^{-1}, 
\y\rightarrow \y^{-1}$ to evaluate the integral (of the same rational 
function) for $|\x|, |\y| < 1$.

We also need to show that the residues of the intermediate integrals at 
zero and infinity are both zero. As we evaluate the first integral, we 
obtain the sum of residues at poles, each of which has $z_1$ equal to a 
monomial in the other variables, possibly with fractional powers. (If 
fractional powers bother you, simply replace each variable by a power of 
itself.) The corresponding residue is obtained by plugging this monomial 
in for $z_1$ in all of the other factors, along with the remaining portion 
of the factor which the pole was extracted out of. Each of these residues 
has denominator a product of terms of the form 1 minus a monomial of the 
variables. As we evaluate each integral, we perform further substitutions 
of monomials in the $x$'s, $y$'s, and later $z$'s for each $z$-variable in 
succession. 

As $z_i$ goes to zero or infinity, each factor in each residue goes to 1, 
0, or a constant depending on whether the exponent of $z_i$ is positive, 
negative, or zero, as in the first integral. We need to show that in each 
residue, one of the factors goes to 0. As $z_i$ goes to infinity, the 
factor $\frac{1}{1-y_i z_i}$, which is unblemished (since it has no other 
$z$-variables), goes to zero. The analysis when $z_i$ goes to 0 is a bit 
trickier. For each previous $z$-variable, representing an equality, we 
have picked one of the factors, corresponding to one $x$ or $y$-variable. 
Substituting for that variable amounts to solving the given equality to 
express that variable in terms of the other variables, and plugging the 
expression into the other equalities to create a new system of equations. 
The residue will be given by the Euler-type generating function of this 
new system of equalities. 

With this formulation, it is clear that if we have a solution 
$\begin{pmatrix} \textbf{x} \\ \textbf{y} \end{pmatrix}$, then the same 
vector with the deleted variable removed will be a solution of the new 
system of equalities; see Example~\ref{computation}.
Taking a vector with all positive entries and 
following it through the process of eliminating the first $i-1$ rows by 
deleting the corresponding columns, we find that the system of equalities 
at the $i$-th step has an all-positive solution. From this, it follows 
that one of the entries in each row must be negative, and in particular 
there must be a negative number in the $i$-th row. The corresponding 
factor of the denominator will have $z_i$ appearing with negative 
exponent, and thus will go to 0 as $z_i$ goes to 0. So each summand has 
residue 0 as $z_i$ goes to 0, and thus the entire $i$-th integrand does.

Therefore, for each integral, we have complementary residues counted in integrals \eqref{firstint} and
\eqref{secondint}, introducing a minus sign; after factoring this out, the two integrals produce
identical rational functions which move on to the next integral. Since we have $m$ integrals, we obtain
that the integrals in \eqref{firstint} and \eqref{secondint} differ by a factor of $(-1)^m$ as desired.
This completes the proof.

It is worth noting that we made the decision to replace the 
$\frac{1}{1-x_i\z^{\a_i}}$ factors in $\theta_{\K}(\x,\y)$ by factors of 
$\frac{x_i\z^{\a_i}}{1-x_i\z^{\a_i}}$ in $\theta_{\K^\circ}(\x,\y)$. Since 
the  facets correspond to the $y$-variables, we do not need strict 
inequalities on the $x$-coordinates in order to compute the generating 
function of $\K^\circ$. Instead, we chose to use this expansion, since it 
yields the correct complementarity statement regarding the poles.

\section{An illustrative example}

In this section, we give an example illustrating the proof of the previous section.

\begin{ex}\label{computation}
Let $\K$ be the pyramidal cone in the positive orthant given by the inequalities:
\begin{eqnarray*}
x_1 - x_2 & \le & 0 \\
x_1 - x_3 & \le & 0 \\
x_2 - 2x_1 & \le & 0 \\
x_3 - 2x_1 & \le & 0.
\end{eqnarray*}
This is a cone with vertex $(0, 0, 0)$ over a square in the $x_1 = 1$ plane with vertices $(1, i, j)$ for $i, j\in \{1, 2\}$.
\end{ex}
The matrix A is then
$\begin{pmatrix} 
1 & -1 & 0 \\
1 & 0 & -1 \\
-2 & 1 & 0 \\
-2 & 0 & 1
\end{pmatrix}$; the modified matrix $A' := (A|I)$ is 
\[
\begin{pmatrix} 
1 & -1 & 0 & 1 & 0 & 0 & 0 \\
1 & 0 & -1 & 0 & 1 & 0 & 0 \\
-2 & 1 & 0 & 0 & 0 & 1 & 0 \\
-2 & 0 & 1 & 0 & 0 & 0 & 1
\end{pmatrix},\]
and the original integrand (to be integrated over $z_1, \ldots, z_4$) is:

\[
\frac{1}{(1-x_1z_1z_2z_3^{-2}z_4^{-2})(1-x_2z_1^{-1}z_3)(1-x_3z_2^{-1}z_4)(1-y_1z_1)(1-y_2z_2)(1-y_3z_3)(1-y_4z_4)}.
\]

We illustrate what happens when we compute a residue in an intermediate integral. 
Suppose we are in the process of computing the residue with respect to $z_1$ at the 
pole corresponding to the first factor. This amounts to substituting $z_1 = 
x_1^{-1}z_2^{-1}z_3^2 z_4^2$ into all of the other factors. Consider the second 
factor, $\left( 1-x_2z_1^{-1}z_3\right)$. This becomes $\left( 1-x_2\left( x_1 z_2 
z_3^{-2}z_4^{-2}\right) z_3\right)$; 
a moment's reflection will reveal that this is equivalent to adding an appropriate 
multiple of the first column of $A'$ (here, 1) to the second column to cancel out 
the element in the first row.

So, in effect, this residue is the Euler-type generating function of the matrix obtained by doing Gaussian elimination, adding a multiple of the first 
column (since we picked the $x_1$ term) to all other columns to eliminate the first row (since we are integrating over $z_1$), and then deleting the 
first column. The only difference is that there are a few $x_1$'s thrown into each term, which are irrelevant for our purposes; recall that the point 
of this statement in the proof was merely to show that in this residue, some term exists with a negative power of each $z_i$.
It suffices to show that this new matrix has a negative entry in each row.

But, as noted in the proof, if we pick a positive solution to the original equation $A'\begin{pmatrix} \textbf{x} \\ \textbf{y} \end{pmatrix} = 0$, 
say $(x_1, x_2, x_3, y_1, y_2, y_3, y_4) = (2, 3, 3, 1, 1, 1, 1)$, then if we remove $x_1$, the new tuple $(x_2, x_3, y_1, y_2, y_3, y_4) = (3, 3, 1, 
1, 1, 1)$ will 
be a solution of the new matrix, which in this case is
\[
\begin{pmatrix} 
1 & -1 & -1 & 1 & 0 & 0 \\
-1 & 0 & 2 & 0 & 1 & 0 \\
-2 & 1 & 2 & 0 & 0 & 1
\end{pmatrix}.\] That the elimination procedure has this property is easily verified; it comprises using the eliminated equation to solve 
for the 
eliminated variable in terms of the other variables, then plugging that into the other equations.

Since there is an all-positive solution to the new set of equations, there must be at least one negative entry in each row.

\section{Extensions and applications}

An interesting extension of Theorem \ref{stanrec}
is the following: Let
\[
  \K_1 = \left\{ \r \in \R_{ \ge 0 }^d : \ \A \, \r \le 0 \text{ and } \B \, \r < 0 \right\}
\]
and
\[
  \K_2 = \left\{ \r \in \R_{ > 0 }^d : \ \A \, \r < 0 \text{ and } \B \, \r \le 0 \right\} ,
\]
that is, $\K_1$ and $\K_2$ are \emph{half-open} cones whose constraints are defined by the same matrix, however, 
those facets (codimension-1 faces) which are contained in $\K_1$ are not 
in $\K_2$ and vice versa. Since by assumption $(0,\ldots,0)$ is the only 
intersection point of the fully closed cone with the coordinate 
hyperplanes, at most one of $\K_1$ and $\K_2$ intersects these 
hyperplanes, so the domains are accurate.


\begin{theorem}\label{stanrecext}
As rational functions, if the set of facets closed in $\K_1$ is contractible, 
$\sigma_{\K_1} \left( 1/\x \right) = (-1)^d \,
\sigma_{ \K_2 } (\x)$.
\end{theorem}

Stanley's original proof~\cite{stanleyreciprocity} showed that this was true 
whenever the set of facets consists of those visible from a certain point outside 
the cone; shelling-based proofs of Theorem~\ref{stanrec} prove this theorem 
whenever the set of closed facets is contractible.

Our proof technique is easily adjusted to this more 
general setting. 

Suppose $\A \in \Z^{ d
\times m }$ has columns $\a_1, \dots, \a_d$, and $\B \in \Z^{ d \times n
}$ has columns $\b_1, \dots, \b_d$, then Lemma \ref{eulerlemma} gives

\[
  \sigma_{ \K_1 } (\x) = \const_{ \z,\w } \left( \frac 1 { \left( 1 - x_1 \z^{ \a_1 } \w^{ \b_1 } \right) \cdots \left( 1 - x_d \z^{ \a_d } \w^{ \b_d } \right) } \frac{ 1 }{ \left( 1 - z_1 \right) \cdots \left( 1 - z_m \right) } \frac{ w_1 }{ 1-w_1 } \cdots \frac{ w_n }{ 1-w_n } \right)
\]
and
\[
  \sigma_{ \K_2 } (\x) = \const_{ \z,\w } \left( \frac{ x_1 \z^{ \a_1 } }{ 
1 - x_1 \z^{ \a_1 } \w^{ \b_1 } } \cdots \frac{ x_d \z^{ \a_d } \w^{ \b_d } }{ 1 - x_d \z^{ \a_d } } \frac{ z_1 }{ 1-z_1 } \cdots \frac{ z_m }{ 1-z_m } \frac{ 1 }{ \left( 1 - w_1 \right) \cdots \left( 1 - w_n \right) } \right) .
\]
The proof that 
\[
  \theta_{\K_1} (\x,\y) = \int \prod_{ j=1 }^{ d } \frac 1 { 1 - x_j \z^{ \a_j } \w^{ \b_j } } \prod_{ k=1 }^{ m } \frac 1 { 1 - y_k z_k } \prod_{ i=1 }^{ n } \frac{ y_{ m+i } w_i }{ 1 - y_{ m+i } w_i } \frac{ d \z }{ \z } \frac{ d \w }{ \w }
\]
and
\[
  \theta_{\K_2} (\x,\y) = \int \prod_{ j=1 }^{ d } \frac{ x_j \z^{ \a_j } \w^{ \b_j } }{ 1 - x_j \z^{ \a_j } \w^{ \b_j } } \prod_{ k=1 }^{ m } \frac{ y_k z_k }{ 1 - y_k z_k } \prod_{ i=1 }^{ n } \frac{ 1 }{ 1 - y_{ m+i } w_i } \frac{ d \z }{ \z } \frac{ d \w }{ \w } \ ,
\]
both defined for $|\x|, |\y| < 1$, satisfy a Stanley-type reciprocity 
identity proceeds along the exact same lines as our proof of Theorem 
\ref{stanrec}. Note that as per the comment at the end of the previous 
section, 
even though the original cones $\K_1$ and $\K_2$ are partially open, we can arbitrarily choose to ``invert'' all of 
the $x$-variables in $\K_2$ and none in $\K_1$ (since $\K_2$ contains no points with any $x$-coordinate equal to 
zero.)

However, this proof cannot work in all cases, since the generalized reciprocity theorem is not true in all cases. What can go wrong is that in the 
partially eliminated matrix, the term (or, rather, all of the terms) with positive $w_i$-exponent can be those of the form $\frac{y_{m+i} 
w_i}{1-y_{m+i}w_i}$, which does not in fact go to zero as $w_i$ goes to infinity. It would be interesting to come up with an elegant charcterization 
of when this happens, which would provide an elegant proof of the generalized Stanley theorem for a subset of half-open cones.

A particular nice application of Theorem \ref{stanrec} concerns the counting function $L_\P (t) := \# \left( t \P \cap \Z^d \right)$ for a rational 
convex polytope $\P$, that is, the convex hull of finitely many points in $\Q^d$. Ehrhart proved in \cite{ehrhartpolynomial} the fundamental 
structural result about $L_\P (t)$, namely that it is a quasi-polynomial in $t$ (for a definition and nice discussion of quasi-polynomials see 
\cite[Chapter 4]{stanleyec1}). Ehrhart conjectured and partially proved the following reciprocity theorem, which was proved by Macdonald 
\cite{macdonald}.

\begin{theorem}[Ehrhart-Macdonald]\label{emrec}
The quasi-polynomials $L_\P$ and $L_{\P^\circ}$ satisfy
$$
  L_\P (-t) = (-1)^{\dim \P} L_{\P^\circ} (t) \ .
$$
\end{theorem} 

\section{Concluding remarks}

1. As already mentioned, most proofs of Theorems \ref{stanrec} and
\ref{emrec} use shellings of a polyhedron or finite additive measures
(see, e.g., \cite{ehrhartbook,macdonald,mcmullenreciprocity}). The only
exceptions we are aware of are proofs via complex analysis (e.g., as above
or in \cite{stanleyreciprocity}) and commutative algebra (see, e.g.,
\cite[Section I.8]{stanleycombcommalg}), as well as a recent proof 
\cite{irrational} by the first author and Frank Sottile using irrational 
decomposition.

2.  It is a fun exercise to deduce Theorem \ref{emrec} from Theorem \ref{stanrec}, for example by considering the
generating function of the $(d+1)$-cone generated by $\left( \v_1 , 1 \right) , \dots, \left( \v_n , 1 \right)$,
where $\v_1, \dots, \v_n$ are the vertices of $\P$, applying Stanley reciprocity, and then specializing the
rational generating functions by setting the first $d$ variables to $1$.

3.  Finally, there exists an extension of Theorem \ref{emrec}
corresponding to Theorem \ref{stanrecext}: one includes some of the facets
of the polytope on one side, and the complementary set of facets on the
other side.



\def\cprime{$'$} \def\cprime{$'$}
\providecommand{\bysame}{\leavevmode\hbox to3em{\hrulefill}\thinspace}
\providecommand{\MR}{\relax\ifhmode\unskip\space\fi MR }
\providecommand{\MRhref}[2]{%
  \href{http://www.ams.org/mathscinet-getitem?mr=#1}{#2}
}
\providecommand{\href}[2]{#2}

\setlength{\parskip}{0cm}

\end{document}